\documentclass{article}

\usepackage[nonatbib,final]{nips_2016}

\usepackage[utf8]{inputenc} 
\usepackage[T1]{fontenc}    
\usepackage{url}            
\usepackage{booktabs}       
\usepackage{amsmath}
\usepackage{amsfonts}       
\usepackage{nicefrac}       
\usepackage{microtype}      
\usepackage{amsthm}
\usepackage{multirow}
\usepackage{color}
\usepackage{soul}
\usepackage{graphicx}
\usepackage{mathtools}
\usepackage{float}
\usepackage{caption}
\usepackage{subcaption}

\usepackage{xcolor}

\usepackage[style=trad-abbrv,backend=bibtex,sorting=none,sortcites=true]{biblatex}
\DeclareFieldFormat[article,inbook,incollection,inproceedings,patent,thesis,unpublished]{title}{{#1}}

\mathtoolsset{centercolon}

\def\dd{\,\mathrm{d}}

\def\Var{\mathrm{Var}}

\newcommand\numberthis{\addtocounter{equation}{1}\tag{\theequation}}
\newtheorem{proposition}{Proposition}
\newtheorem{theorem}{Theorem}
\newtheorem*{theorem*}{Theorem}

\renewcommand\arraystretch{1.2}

\title{Probabilistic Linear Multistep Methods}

\author{
  Onur Teymur\\
  Department of Mathematics\\
  Imperial College London\\
  \texttt{o@teymur.uk} \\
  \And
  Konstantinos Zygalakis \\
  School of Mathematics  \\
  University of Edinburgh \\
  \texttt{k.zygalakis@ed.ac.uk} \\
  \And
  Ben Calderhead \\
  Department of Mathematics \\
  Imperial College London \\
  \texttt{b.calderhead@imperial.ac.uk} \\
}

\begin{document}
\maketitle

\begin{abstract}
	We present a derivation and theoretical investigation of the Adams-Bashforth and Adams-Moulton family of linear multistep methods for solving ordinary differential equations, starting from a Gaussian process (GP) framework. In the limit, this formulation coincides with the classical deterministic methods, which have been used as higher-order initial value problem solvers for over a century. Furthermore, the natural probabilistic framework provided by the GP formulation allows us to derive \textit{probabilistic} versions of these methods, in the spirit of a number of other probabilistic ODE solvers presented in the recent literature \mbox{\cite{chkrebtii16,hennig14,schober14,conrad16}}. In contrast to higher-order Runge-Kutta methods, which require multiple intermediate function evaluations per step, Adams family methods make use of previous function evaluations, so that increased accuracy arising from a higher-order multistep approach comes at very little additional computational cost. We show that through a careful choice of covariance function for the GP, the posterior mean and standard deviation over the numerical solution can be made to exactly coincide with the value given by the deterministic method and its local truncation error respectively.  We provide a rigorous proof of the convergence of these new methods, as well as an empirical investigation (up to fifth order) demonstrating their convergence rates in practice.
\end{abstract}

\section{Introduction}
Numerical solvers for differential equations are essential tools in almost all disciplines of applied mathematics, due to the ubiquity of real-world phenomena described by such equations, and the lack of exact solutions to all but the most trivial examples.
The performance -- speed, accuracy, stability, robustness -- of the numerical solver is of great relevance to the practitioner. This is particularly the case if the computational cost of accurate solutions is significant, either because of high model complexity or because a high number of repeated evaluations are required (which is typical if an ODE model is used as part of a statistical inference procedure, for example). A field of work has emerged which seeks to quantify this performance -- or indeed lack of it -- by modelling the numerical errors \textit{probabilistically}, and thence trace the effect of the chosen numerical solver through the entire computational pipeline \cite{kennedy01}. The aim is to be able to make meaningful quantitative statements about the uncertainty present in the resulting scientific or statistical conclusions.

Recent work in this area has resulted in the development of probabilistic numerical methods, first conceived in a very general way in \cite{diaconis88}. An recent summary of the state of the field is given in \cite{hennig15}. The particular case of ODE solvers was first addressed in \cite{skilling93}, formalised and extended in \cite{chkrebtii16,hennig14,schober14} with a number of theoretical results recently given in \cite{conrad16}. The present paper modifies and extends the constructions in \cite{chkrebtii16,conrad16} to the multistep case, improving the order of convergence of the method but avoiding the simplifying linearisation of the model required by the approaches of \cite{hennig14,schober14}. Furthermore we offer extensions to the convergence results in \cite{conrad16} to our proposed method and give empirical results confirming convergence rates which point to the practical usefulness of our higher-order approach without significantly increasing computational cost.

\subsection {Mathematical setup}

We consider an Initial Value Problem (IVP) defined by an ODE
\begin{equation}\label{ode defn}
\frac{\dd}{\dd t} y(t,\theta) = f(y(t,\theta),t), \qquad \qquad y(t_0,\theta) = y_0
\end{equation}
Here $y(\cdot,\theta): \mathbb{R}^+ \rightarrow \mathbb{R}^d$ is the solution function, $f:\mathbb{R}^d \times \mathbb{R}^+ \rightarrow \mathbb{R}^d$ is the vector-valued function that defines the ODE, and $y_0 \in \mathbb{R}^d$ is a given vector called the initial value. The dependence of $y$ on an $m$-dimensional parameter $\theta \in \mathbb{R}^m$ will be relevant if the aim is to incorporate the ODE into an inverse problem framework, and this parameter is of scientific interest. Bayesian inference under this setup (see \cite{girolami08}) is covered in most of the other treatments of this topic but is not the main focus of this paper; we therefore suppress $\theta$ for the sake of clarity.

Some technical conditions are required in order to justify the existence and uniqueness of solutions to \eqref{ode defn}. We assume that $f$ is evaluable point-wise given $y$ and $t$ and also that it satisfies the Lipschitz condition in $y$, namely $||f(y_1,t) - f(y_2,t)|| \leq L_f||y_1 - y_2||$ for some $L_f \in \mathbb{R}^+$ and all $t,y_1$ and $y_2$; and also is continuous in $t$. These conditions imply the existence of a unique solution, by a classic result usually known as the Picard-Lindel{\"o}f Theorem \cite{iserles08}.

We consider a finite-dimensional discretisation of the problem, with our aim being to numerically generate an $N$-dimensional vector\footnote{The notation $y_{0:N}$ denotes the vector $\left(y_0,\dots,y_N\right)$, and analogously $t_{0:N}$, $f_{0:N}$ etc.} $y_{1:N}$ approximating the true solution $y(t_{1:N})$ in an appropriate sense. Following \cite{chkrebtii16}, we consider the joint distribution of $y_{1:N}$ and the auxiliary variables $f_{0:N}$ (obtained by evaluating the function $f$), with each $y_i$ obtained by sequentially conditioning on previous evaluations of $f$. A basic requirement is that the marginal mean of $y_{1:N}$ should correspond to some deterministic iterative numerical method operating on the grid $t_{1:N}$. In our case this will be a linear multistep method (LMM) of specified type.\;\footnote{We argue that the connection to some specific deterministic method is a desirable feature, since it aids interpretability and allows much of the well-developed theory of IVP solvers to be inherited by the probabilistic solver. This is a particular strength of the formulation in \cite{conrad16} which was lacking in all previous works.}

Firstly we telescopically factorise the joint distribution as follows:
\begin{equation}
	p(y_{1:N},f_{0:N}|y_0) =
	p(f_0|y_0) \prod_{i=0}^{N-1} p(y_{i+1}|y_{0:i},f_{0:i})\; p(f_{i+1}|y_{0:i+1},f_{0:i}) \label{telescope}
\end{equation}
We can now make simplifying assumptions about the constituent distributions. Firstly since we have assumed that $f$ is evaluable point-wise given $y$ and $t$,
\begin{equation}
p(f_i|y_i,\dots) = p(f_i|y_i) = \delta_{f_i}\big(f(y_i,t_i)\big),
\end{equation}
 which is a Dirac-delta measure equivalent to simply performing this evaluation deterministically. Secondly, we assume a finite moving window of dependence for each new state -- in other words $y_{i+1}$ is only allowed to depend on $y_i$ and $f_i,f_{i-1},\dots,f_{i-(s-1)}$ for some $s \in \mathbb{N}$. This corresponds to the inputs used at each iteration of the $s$-step Adams-Bashforth method. For $i<s$ we will assume dependence on only those derivative evaluations up to $i$; this initialisation detail is discussed briefly in Section \ref{implementation}. Strictly speaking, $f_N$ is superfluous to our requirements (since we already have $y_N$) and thus we can rewrite \eqref{telescope} as
\begin{align}
	p(y_{1:N},f_{0:N-1}|y_0) &= \prod_{i=0}^{N-1} p(f_{i}|y_{i}) \: p (y_{i+1}|y_i,f_{\mathrm{max}(0,i-s+1):i}) \\
	&= \prod_{i=0}^{N-1} \delta_{f_i}(f(y_i,t_i)) \: \underbrace{p(y_{i+1}|y_i,f_{\mathrm{max}(0,i-s+1):i})}_{\ast} \label{product telescope}
\end{align}
The conditional distributions $\ast$ are the primary objects of our study -- we will define them by constructing a particular Gaussian process prior over all variables, then identifying the appropriate (Gaussian) conditional distribution. Note that a simple modification to the decomposition \eqref{telescope} allows the same set-up to generate an $(s+1)$-step Adams-Moulton iterator\footnote{The convention is that the number of steps is equal to the total number of derivative evaluations used in each iteration, hence the $s$-step AB and $(s+1)$-step AM methods both go `equally far back'.} -- the implicit multistep method where $y_{i+1}$ depends in addition on $f_{i+1}$. At various stages of this paper this extension is noted but omitted for reasons of space -- the collected results are given in Appendix \ref{am extension}.

\subsection*{Linear multistep methods}

We give a very short summary of Adams family LMMs and their conventional derivation via interpolating polynomials. For a fuller treatment of this well-studied topic we refer the reader to the comprehensive references \cite{iserles08,hairer08,butcher08}. Using the usual notation we write $y_i$ for the numerical estimate of the true solution $y(t_i)$, and $f_i$ for the estimate of $f(t_i) \equiv y'(t_i)$.

The classic $s$-step Adams-Bashforth method calculates $y_{i+1}$ by constructing the unique polynomial $P_i(\omega) \in \mathbb{P}_{s-1}$ interpolating the points $\{f_{i-j}\}_{j=0}^{s-1}$. This is given by \mbox{Lagrange's} method as
\begin{equation}
P_i(\omega) = \sum_{j=0}^{s-1} \ell_j^{\,0:s-1}(\omega) f_{i-j}  \qquad\qquad \ell_j^{\,0:s-1}(\omega) = \prod_{\substack{k=0 \\ k \neq j}}^{s-1} \frac{\omega - t_{i-k}}{t_{i-j}-t_{i-k}} \label{lagrange interpolation}
\end{equation}
The $\ell_j^{\,0:s-1}(\omega)$ are known as Lagrange polynomials, have the property that $\ell_p^{\,0:s-1}(t_{i-q}) = \delta_{pq}$, and form a basis for the space $\mathbb{P}_{s-1}$ known as the Lagrange basis. The Adams-Bashforth iteration then proceeds by writing the integral version of \eqref{ode defn} as $ y(t_{i+1}) - y(t_{i}) \equiv \int_{t_i}^{t_{i+1}} f(y,t) \dd t $ and approximating the function under the integral by the extrapolated interpolating polynomial to give
\begin{equation}
y_{i+1} - y_{i} \approx \int_{t_i}^{t_{i+1}} P_i(\omega) \dd\omega = h \sum_{j=0}^{s-1} \beta^{AB}_{j,s} f_{i-j} \label{ab def}
\end{equation}
where $h = t_{i+1} - t_i$ and the $\beta^{AB}_{j,s} \equiv h^{-1}\int_0^h \ell_j^{\,0:s-1}(\omega)\dd\omega$ are the Adams-Bashforth coefficients for order $s$, all independent of $h$ and summing to 1. Note that if $f$ is a polynomial of degree $s-1$ (so $y(t)$ is a polynomial of degree $s$) this procedure will give the next solution value \textit{exactly}. Otherwise the extrapolation error in $f_{i+1}$ is of order $O(h^s)$ and in $y_{i+1}$ (after an integration) is of order $O(h^{s+1})$. So the local truncation error is $O(h^{s+1})$ and the global error $O(h^{s})$ \cite{iserles08}.

Adams-Moulton methods are similar except that the polynomial $Q_i(\omega) \in \mathbb{P}_{s}$ interpolates the $s+1$ points $\{f_{i-j}\}_{j=-1}^{s-1}$. The resulting equation analogous to \eqref{ab def} is thus an implicit one, with the unknown $y_{i+1}$ appearing on both sides. Typically AM methods are used in conjunction with an AB method of one order lower, in a `predictor-corrector' arrangement. Here, a predictor value ${y}^\ast_{i+1}$ is calculated using an AB step; this is then used to estimate $f^\ast_{i+1} = f(y^\ast_{i+1})$; and finally an AM step uses this value to calculate $y_{i+1}$. We again refer the reader to Appendix \ref{am extension} for details of the AM construction.

\section{Derivation of Adams family LMMs via Gaussian processes}
We now consider a formulation of the Adams-Bashforth family starting from a Gaussian process framework and then present a probabilistic extension. We fix a joint Gaussian process prior over $y_{i+1},y_{i},f_i,f_{i-1},\dots,f_{i-s+1}$ as follows. We define two vectors of functions $\phi(\omega)$ and $\Phi(\omega)$ in terms of the Lagrange polynomials $\ell^{\,0:s-1}_j(\omega)$ defined in \eqref{lagrange interpolation} as
\begin{align}
\phi(\omega) &= \begin{pmatrix} 0 & \ell_0^{\,0:s-1}(\omega) & \ell_1^{\,0:s-1}(\omega) &  \vphantom{\displaystyle\int}\dots & \ell_{s-1}^{\,0:s-1}(\omega) \end{pmatrix}^T \label{phi}\\
\Phi(\omega) &= \int \phi(\omega) \dd \omega  = \begin{pmatrix} 1 & \displaystyle\int\ell_0^{\,0:s-1}(\omega)\dd\omega\; & \dots &\displaystyle\int\ell_{s-1}^{\,0:s-1}(\omega)\dd\omega \end{pmatrix}^T \label{Phi}
\end{align}
The elements (excluding the first) of $\phi(\omega)$ form a basis for $\mathbb{P}_{s-1}$ and the elements of $\Phi(\omega)$ form a basis for $\mathbb{P}_{s}$. The initial $0$ in $\phi(\omega)$ is necessary to make the dimensions of the two vectors equal, so we can correctly define products such as $\Phi(\omega)^T\phi(\omega)$ which will be required later. The first element of $\Phi(\omega)$ can be any non-zero constant $C$; the analysis later is unchanged and we therefore take $C = 1$.

Since we will solely be interested in values of the argument $\omega$ corresponding to discrete equispaced time-steps $t_{j} - t_{j-1} = h$ indexed relative to the current time-point $t_i = 0$, we will make our notation more concise by writing $\phi_{i+k}$ for $\phi(t_{i+k})$, and similarly $\Phi_{i+k}$ for $\Phi(t_{i+k})$. We now use these vectors of basis functions to define a joint Gaussian process prior as follows:
\newpage
\begin{equation}
\hspace{-0mm}
\begin{pmatrix}
	y_{i+1}\\ y_{i}\\ f_i\\ f_{i-1}\\ \vdots \\ f_{i-s+1}
\end{pmatrix} = \mathcal{N} \left[ \begin{pmatrix}
	0 \\ 0 \\ 0 \\ 0 \\ \vdots \\ 0 
\end{pmatrix}, \begin{pmatrix}
	\Phi_{i+1}^T\Phi_{i+1} & \Phi_{i+1}^T\Phi_{i} & \Phi_{i+1}^T\phi_{i} & \cdots & \Phi_{i+1}^T\phi_{i-s+1}\\
	
	\Phi_{i}^T\Phi_{i+1} & \Phi_{i}^T\Phi_{i} & \Phi_{i}^T\phi_{i} & \cdots & \Phi_{i}^T\phi_{i-s+1}\\
	
	\phi_{i}^T\Phi_{i+1} & \phi_{i}^T\Phi_{i} &\phi_i^T\phi_i & \dots & \phi_i^T\phi_{i-s+1} &  \\
	
	\phi_{i-1}^T\Phi_{i+1} & \phi_{i-1}^T\Phi_{i}& \phi^T_{i-1}\phi_i & \dots & \phi_{i-1}^T\phi_{i-s+1} \\
	
	\vdots & \vdots & \vdots & \ddots &\vdots\\
	\phi_{i-s+1}^T\Phi_{i+1}&
	\phi_{i-s+1}^T\Phi_{i} & \phi_{i-s+1}^T\phi_{i} & \dots & \phi_{i-s+1}^T\phi_{i-s+1}
\end{pmatrix} \right] \label{ab gp}
\end{equation}

This construction works because $y' = f$ and differentiation is a linear operator; the rules for the transformation of the covariance elements is given in Section 9.4 of \cite{rasmussen06} and can easily be seen to correspond to the defined relationship between $\phi(\omega)$ and $\Phi(\omega)$.

Recalling the decomposition in \eqref{product telescope}, we are interested in the conditional distribution $p(y_{i+1}|y_i,f_{i-s+1:i})$. This is also Gaussian, with mean and covariance given by the standard formulae for Gaussian conditioning. This construction now allows us to state the following result:

\begin{proposition}
The conditional distribution $p(y_{i+1}|y_i,f_{i-s+1:i})$ under the Gaussian process prior given in \eqref{ab gp}, with covariance kernel basis functions as in \eqref{phi} and \eqref{Phi}, is a $\delta$-measure concentrated on the $s$-step Adams-Bashforth predictor $y_i + h \sum_{j=0}^{s-1} \beta^{AB}_{j,s} f_{i-j}$. \label{ab det theorem}
\end{proposition}
The proof of this proposition is given in Appendix \ref{p1 proof}.

Because of the natural probabilistic structure provided by the Gaussian process framework, we can augment the basis function vectors $\phi(\omega)$ and $\Phi(\omega)$ to generate a conditional distribution for $y_{i+1}$ that has non-zero variance. By choosing a particular form for this augmented basis we can obtain an expression for the standard deviation of $y_{i+1}$ that is exactly equal to the leading-order local truncation error of the corresponding deterministic method.

We will expand the vectors $\phi(\omega)$ and $\Phi(\omega)$ by one component, chosen so that the new vector comprises elements that span a polynomial space of order one greater than before. Define the augmented bases $\phi^+(\omega)$ and $\Phi^+(\omega)$ as
\begin{align}
\phi(\omega)^+ &= \begin{pmatrix} \phantom{\displaystyle\int} 0 & \ell_0^{\,0:s-1}(\omega) & \ell_1^{\,0:s-1}(\omega) & \dots & \ell_{s-1}^{\,0:s-1}(\omega) & \alpha h^{s}\ell_{-1}^{\,-1:s-1}(\omega)\end{pmatrix}^T \label{phiplus}\\
\Phi(\omega)^+ &= \begin{pmatrix} 1 & \displaystyle\int\ell_0^{\,0:s-1}(\omega)\dd\omega\; & \dots &\displaystyle\int\ell_{s-1}^{\,0:s-1}(\omega)\dd\omega & \displaystyle \int \alpha h^{s}\ell_{-1}^{\,-1:s-1}(\omega) \dd\omega\end{pmatrix}^T \label{Phiplus}
\end{align}
The additional term at the end of $\phi^+(\omega)$ is the polynomial of order $s$ which arises from interpolating $f$ at $s+1$ points (with the additional point at $t_{i+1}$) and choosing the basis function corresponding to the root at $t_{i+1}$, scaled by $\alpha h^{s}$ with $\alpha$ a positive constant whose role will be explained in the next section. The elements of these vectors span $\mathbb{P}_{s}$ and $\mathbb{P}_{s+1}$ respectively. With this new basis we can give the following result:

\begin{proposition}
The conditional distribution $p(y_{i+1}|y_i,f_{i-s+1:i})$ under the Gaussian process prior given in \eqref{ab gp}, with covariance kernel basis functions as in \eqref{phiplus} and \eqref{Phiplus}, is Gaussian with mean equal to the $s$-step Adams-Bashforth predictor $y_i~+~h \sum_{j=0}^{s-1} \beta^{AB}_{j,s} f_{i-j}$ and, setting $\alpha = y^{(s+1)}(\eta)$ for some $\eta \in (t_{i-s+1},t_{i+1})$, standard deviation equal to its local truncation error.
\end{proposition}
The proof is given in Appendix \ref{p2 proof}. In order to de-mystify the construction, we now exhibit a concrete example for the case $s=3$. The conditional distribution of interest is $p(y_{i+1}|y_i,f_i,f_{i-1},f_{i-2}) \equiv p(y_{i+1}|y_i,f_{i:i-2})$. In the deterministic case, the vectors of basis functions become
\begin{align*}
\phi(\omega)_{s=3} &= \begin{pmatrix} 0 & \dfrac{(\omega+h)(\omega+2h)}{2h^2} & \dfrac{\omega(\omega+2h)}{-h^2} & \dfrac{\omega(\omega+h)}{2h^2} \end{pmatrix} \\
\Phi(\omega)_{s=3} &= \begin{pmatrix} 1 & \dfrac{\omega \left( 2\omega^2 + 9h\omega + h^2 \right)}{12h^2} & \dfrac{\omega^2 \left( \omega + 3h \right)}{-3h^2} & \dfrac{\omega^2 \left( 2\omega + 3h \right)}{12h^2} \end{pmatrix}
\end{align*}
\newpage
and simple calculations give that 
\begin{align*}
\mathbb{E}(y_{i+1}|y_i,f_{i:i-2}) = y_i + h\left( \dfrac{23}{12}f_i -\dfrac{4}{3}f_{i-1} + \dfrac{5}{12}f_{i-2} \right) \qquad \quad
\mathrm{\Var}(y_{i+1}|y_i,f_{i:i-2}) = 0
\end{align*}
The probabilistic version follows by setting
\begin{align*}
\phi^+(\omega)_{s=3} &= \begin{pmatrix} 0 & \dfrac{(\omega+h)(\omega+2h)}{2h^2} & \dfrac{\omega(\omega+2h)}{-h^2} & \dfrac{\omega(\omega+h)}{2h^2} & \dfrac{\alpha\omega(\omega+h)(\omega+2h)}{6} \end{pmatrix} \\
\Phi^+(\omega)_{s=3} &= \begin{pmatrix} 1 &\dfrac{\omega \left( 2\omega^2 + 9h\omega + h^2 \right)}{12h^2} & \dfrac{\omega^2 \left( x + 3h \right)}{-3h^2} & \dfrac{\omega^2 \left( 2\omega + 3h \right)}{12h^2} & \dfrac{\alpha\omega^2 (\omega + 2h)^2}{24} \end{pmatrix}
\end{align*}
and further calculation shows that
\begin{align*}
	\mathbb{E}(y_{i+1}|y_i,f_{i:i-2}) = y_i + h\left( \dfrac{23}{12}f_i -\dfrac{4}{3}f_{i-1} + \dfrac{5}{12}f_{i-2} \right) \qquad
\mathrm{\Var}(y_{i+1}|y_i,f_{i:i-2}) = \left( \dfrac{3 h^4 \alpha}{8} \right)^2
\end{align*}

An entirely analogous argument can be shown to reproduce and probabilistically extend the implicit Adams-Moulton scheme. The Gaussian process prior now includes $f_{i+1}$ as an additional variable and the correlation structure and vectors of basis functions are modified accordingly. The required modifications are given in Appendix \ref{am extension} and a explicit derivation for the 4-step AM method is given in Appendix \ref{s4}. 

\subsection{The role of $\alpha$} \label{alpha}
Replacing $\alpha$ in \eqref{phiplus} by $y^{(s+1)}(\eta)$, with $\eta \in (t_{i-s+1},t_{i+1})$, makes the variance of the integrator coincide exactly with the local truncation error of the underlying deterministic method.\footnote{We do not claim that this is the only possible way of modelling the numerical error in the solver. The question of how to do this accurately is an open problem in general, and is particularly challenging in the multi-dimensional case. In many real world problems different noise scales will be appropriate for different dimensions and -- especially in `hierarchical' models arising from higher-order ODEs -- non-Gaussian noise is to be expected. That said, the Gaussian assumption as a first order approximation for numerical error is present in virtually all work on this subject and goes all the way back to \cite{skilling93}. We adopt this premise throughout, whilst noting this interesting unresolved issue.}
\iffalse Our approach is heuristically justifiable due to its immediate connection with the local truncation error of the deterministic method.} \fi

This is of course of limited utility unless higher derivatives of $y(t)$ are available, and even if they are, $\eta$ is itself unknowable in general. However it is possible to estimate the integrator variance in a systematic way by using backward difference approximations \cite{fornberg88} to the required derivative at $t_{i+1}$. We show this by expanding the $s$-step Adams-Bashforth iterator as
\begin{align*}
y_{i+1} &= y_i + h\textstyle\sum_{j=0}^{s-1}\beta^{AB}_{j,s}f_{i-j} + h^{s+1}C_s^{AB}y^{(s+1)}(\eta) \qquad \qquad \qquad \qquad \eta \in [t_{i-s+1},t_{i+1}] \\
&= y_i + h\textstyle\sum_{j=0}^{s-1}\beta^{AB}_{j,s}f_{i-j} + h^{s+1}C_s^{AB}y^{(s+1)}(t_{i+1}) + O(h^{s+2})\\
&= y_i + h\textstyle\sum_{j=0}^{s-1}\beta^{AB}_{j,s}f_{i-j} + h^{s+1}C_s^{AB}f^{(s)}(t_{i+1}) + O(h^{s+2}) \qquad \quad \text{since } y' = f\\
&= y_i + h\textstyle\sum_{j=0}^{s-1}\beta^{AB}_{j,s}f_{i-j} + h^{s+1}C_s^{AB}\left[h^{-s}\textstyle\sum_{k=0}^{s-1+p}\delta_{k,s-1+p}f_{i-k} + O(h^p) \right] + O(h^{s+2})\\
&= y_i + h\textstyle\sum_{j=0}^{s-1}\beta^{AB}_{j,s}f_{i-j} + h C_s^{AB}\textstyle\sum_{k=0}^{s}\delta_{k,s}f_{i-k} + O(h^{s+2}) \qquad \;\text{if we set } p = 1 \numberthis\label{bdc}
\end{align*}
where $\beta^{AB}_{\cdot,s}$ are the set of coefficients and $C_s^{AB}$ the local truncation error constant for the $s$-step Adams-Bashforth method, and $\delta_{\cdot,s-1+p}$ are the set of backward difference coefficients for estimating the $s$th derivative of $f$ to order $O(h^p)$ \cite{fornberg88}.

In other words, the constant $\alpha$ can be substituted with $h^{-s}\sum_{k=0}^{s}\delta_{k,s}f_{i-k}$, using already available function values and to adequate order. It is worth noting that collecting the coefficients $\beta^{AB}_{\cdot,s}$ and $\delta_{\cdot,s}$ results in an expression equivalent to the Adams-Bashforth method of order $s+1$ and therefore, this procedure is in effect employing two integrators of different orders and estimating the truncation error from the difference of the two.\footnote{An explicit derivation of this for $s=3$ is given in Appendix \ref{bdc expansion}.} This principle is similar to the classical Milne Device \cite{butcher08}, which pairs an AB and and AM iterator to achieve the same thing. Using the Milne Device to generate a value for the error variance is also straightforward within our framework, but requires two evaluations of $f$ at each iteration (one of which immediately goes to waste) instead of the approach presented here, which only requires one.
\section{Convergence of the probabilistic Adams-Bashforth integrator}
We now give the main result of our paper, which demonstrates that the convergence properties of the probabilistic Adams-Bashforth integrator match those of its deterministic counterpart.
\setcounter{theorem}{2}
\begin{theorem}
Consider the $s$-step deterministic Adams-Bashforth integrator given in Proposition 1, which is of order $s$. Then the probabilistic integrator constructed in Proposition 2 has the same mean square error as its deterministic counterpart. In particular 
\[
\max_{0 \leq kh \leq T}\mathbb{E}|Y_{k}-y_{k}|^2 \leq K h^{2s}
\]
where $Y_k \equiv y(t_k)$ denotes the true solution, $y_k$ the numerical solution, and K is a positive real number depending on $T$ but independent of $h$.
\end{theorem}
The proof of this theorem is given in Appendix \ref{t3 proof}, and follows a similar line of reasoning to that given for a one-step probabilistic Euler integrator in \cite{conrad16}. In particular, we deduce the convergence of the 
algorithm by extrapolating from the local error. The additional 
complexity arises due to the presence 
of the stochastic part, which means we cannot rely directly on the 
theory of difference equations and the representations of 
their solutions. Instead, following \cite{buckwar06}, we rewrite the defining $s$-step recurrence equation as a one-step recurrence equation in a higher dimensional space.

\section{Implementation} \label{implementation}
We now have an implementable algorithm for an $s$-step probabilistic Adams-Bashforth integrator. Firstly, an accurate initialisation is required for the first $s$ iterations -- this can be achieved with, for example, a Runge-Kutta method of sufficiently high order.\footnote{We use a (packaged) adaptive Runge-Kutta-Fehlberg solver of 7th order with 8th order error control.} Secondly, at iteration $i$, the preceding $s$ stored function evaluations are used to find the posterior mean and variance of $y_{i+1}$. The integrator then advances by generating a realisation of the posterior measure derived in Proposition 2. Following \cite{chkrebtii16}, a Monte Carlo repetition of this procedure with different random seeds can then be used as an effective way of generating propagated uncertainty estimates at any time $0<T<\infty$. 
\subsection{Example -- Chua circuit} \label{chua section}
The Chua circuit \cite{chua92} is the simplest electronic circuit that exhibits chaotic behaviour, and has been the subject of extensive study -- in both the mathematics and electronics communities -- for over 30 years. Readers interested in this rich topic are directed to \cite{chua07} and the references therein. The defining characteristic of chaotic systems is their unpredictable long-term sensitivity to tiny changes in initial conditions, which also manifests itself in the sudden amplification of error introduced by any numerical scheme. It is therefore of interest to understand the limitations of a given numerical method applied to such a problem -- namely the point at which the solution can no longer be taken to be a meaningful approximation of the ground truth. Probabilistic integrators allow us to do this in a natural way \cite{chkrebtii16}.

The Chua system is given by $x' = \alpha(y- (1+h_1) x - h_3 x^3)$, $y' = x-y+z$, $z' = -\beta y - \gamma z$. We use parameter values $\alpha = -1.4157$, $\beta = 0.02944201$, $\gamma = 0.322673579$, $h_1 = -0.0197557699$, $h_3 = -0.0609273571$ and initial conditions $x_0 = 0$, $y_0 = 0.003$, $z_0 = 0.005$. This particular choice is taken from `Attractor CE96' in \cite{bilotta08}.
\begin{figure}[h]
\captionsetup{width=0.9\textwidth,justification=centering}
  \centering
  \includegraphics[width=\textwidth,height=6.5cm]{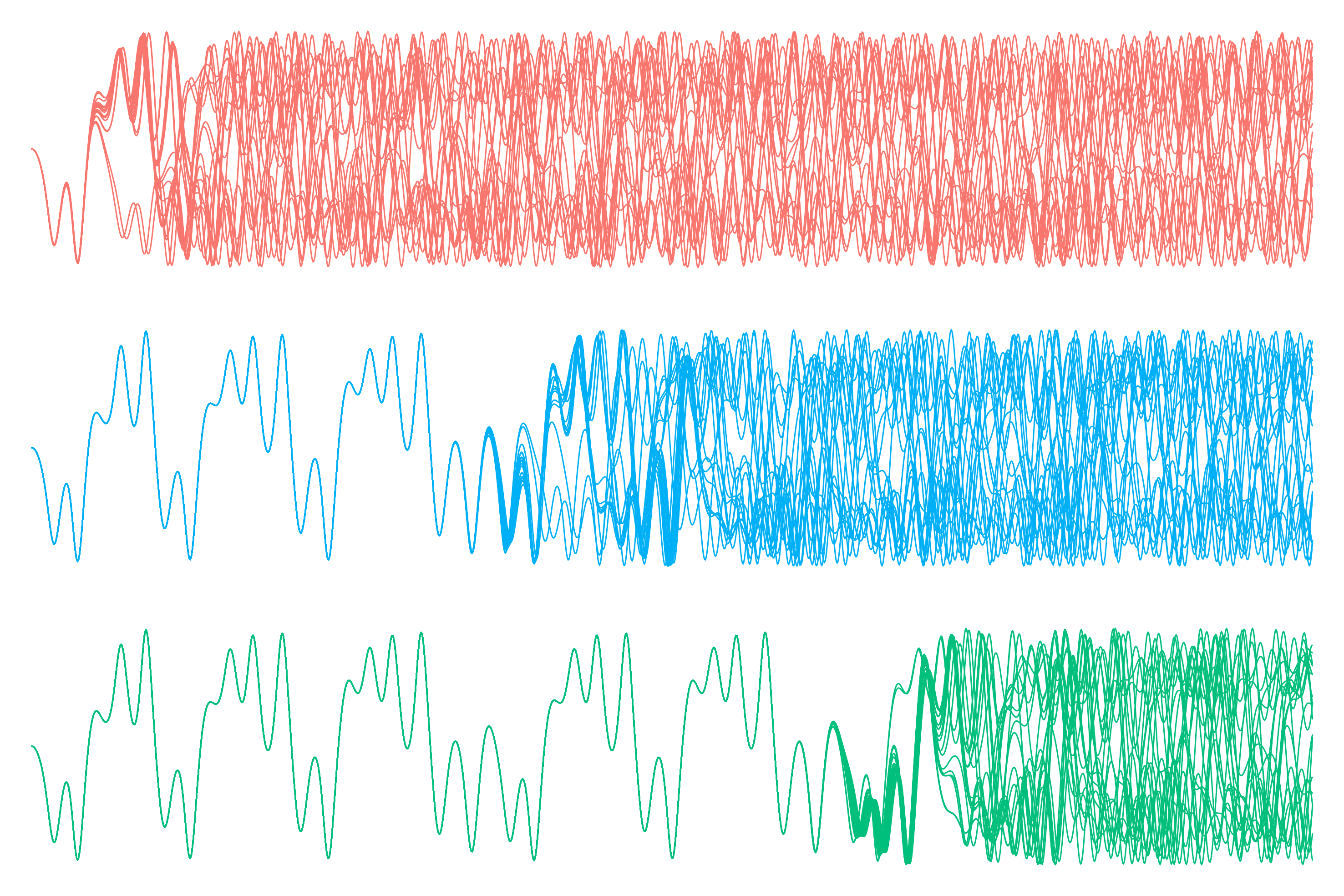}
  \caption{Time series for the $x$-component in the Chua circuit model described in Section \ref{chua section}, solved 20 times for $0\leq t \leq 1000$ using an $s$-step probabilistic AB integrator with $s=1$ (top), $s=3$ (middle), $s=5$ (bottom). Step-size remains $h=0.01$ throughout. Wall-clock time for each simulation was close to constant ($\pm 10$ per cent -- the difference primarily accounted for by the RKF initialisation procedure).} \label{chua plot}
\end{figure}
Using the probabilistic version of the Adams-Bashforth integrator with $s>1$, it is possible to delay the point at which numerical path diverges from the truth, \textit{with effectively no additional evaluations of $f$ required} compared to the one-step method. This is demonstrated in Figure~\ref{chua plot}. Our approach is therefore able to combine the benefits of classical higher-order methods with the additional insight into solution uncertainty provided by a probabilistic method.

\subsection{Example -- Lotka-Volterra model} \label{lv section}
We now apply the probabilistic integrator to a simple periodic predator-prey model given by the system $x' = \alpha x - \beta x y$, $y' = \gamma x y - \delta y$ for parameters $\alpha = 1$, $\beta = 0.3$, $\gamma = 1$ and $\delta = 0.7$. We demonstrate the convergence behaviour stated in Theorem 3 empirically.

\begin{figure}[H]
\captionsetup{width=0.9\textwidth,justification=centering}
\hspace{-4mm}  \includegraphics[width=1.05\linewidth]{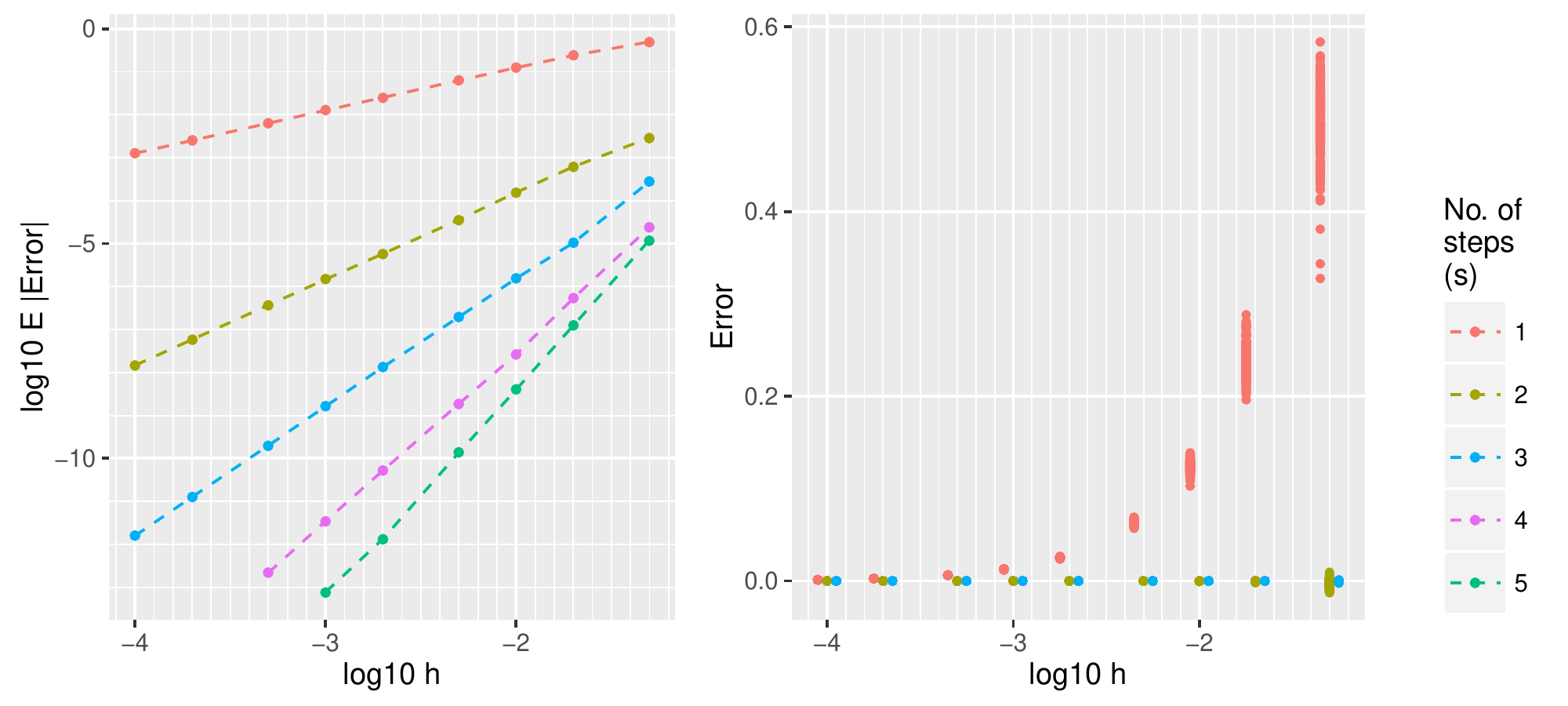}
  \caption{Empirical error analysis for the $x$-component of 200 realisations of the probabilistic AB integrator as applied to the Lotka-Volterra model described in Section \ref{lv section}. The left-hand plot shows the convergence rates for AB integrators of orders 1-5, while the right-hand plot shows the distribution of error around zero for integrators of orders 1-3.}
  \label{lv}
\end{figure}

The left-hand plot in Figure \ref{lv} shows the sample mean of the absolute error of 200 realisations of the probabilistic integrator plotted against step-size, on a log-log scale. The differing orders of convergence of the probabilistic integrators are easily deduced from the slopes of the lines shown.

The right-hand plot shows the actual error value (no logarithm or absolute value taken) of the same 200 realisations, plotted individually against step-size. This plot shows that the error in the one-step integrator is consistently positive, whereas for two- and three-step integrators is approximately centred around 0. (This is also visible with the same data if the plot is zoomed to more closely examine the range with small $h$.) Though this phenomenon can be expected to be somewhat problem-dependent, it is certainly an interesting observation which may have implications for bias reduction in a Bayesian inverse problem setting.

\section{Conclusion}
 We have given a derivation of the Adams-Bashforth and Adams-Moulton families of linear multistep ODE integrators, making use of a Gaussian process framework, which we then extend to develop their probabilistic counterparts.  

We have shown that the derived family of probabilistic integrators result in a posterior mean at each step that exactly coincides with the corresponding deterministic integrator, with the posterior standard deviation equal to the deterministic method's local truncation error. We have given the general forms of the construction of these new integrators to arbitrary order. Furthermore, we have investigated their theoretical properties and provided a rigorous proof of their rates of convergence, Finally we have demonstrated the use and computational efficiency of probabilistic Adams-Bashforth methods by implementing the solvers up to fifth order and providing example solutions of a chaotic system, and well as empirically verifying the convergence rates in a Lotka-Voltera model.

We hope the ideas presented here will add to the arsenal of any practitioner who uses numerical methods in their scientific analyses, and contributes a further tool in the emerging field of probabilistic numerical methods.

\renewcommand*{\bibfont}{\small}
\printbibliography
\vspace{-2mm}
\let\svthefootnote\thefootnote
\let\thefootnote\relax\footnote{KZ was partially supported by a grant from the Simons Foundation. Part of this work was done during the author's stay at the Newton Institute for the programme \textit{Stochastic Dynamical Systems in Biology: Numerical Methods and Applications.}}
\addtocounter{footnote}{-1}\let\thefootnote\svthefootnote
\clearpage

\newpage
\section*{Appendices}
\pagenumbering{roman}
\renewcommand{\thesubsection}{\Alph{subsection}}

\subsection{Proof of Proposition 1} \label{p1 proof}
Recall that $h = t_{j} - t_{j-1}$ for all $j$. Straightforward substitutions into the definitions give that $\phi_i \equiv \phi(0) = (0,1,0,\dots,0)$, $\phi_{i-1} \equiv \phi(-h) = (0,0,1,\dots,0)$ etc. and hence $\phi_{i-p}^T\phi_{i-q} = \delta_{pq}$, for all $0 \leq p,q\leq s-1$. Furthermore $\Phi_i \equiv \Phi(0) = (1,0,0,\dots,0)$ since every component of $\Phi(\omega)$ bar the first is a polynomial of degree $s$ with a factor $\omega$. Finally 
\[
\Phi_{i+1} \equiv \Phi(h) = \begin{pmatrix} 1 & \displaystyle\int^h_0 \ell_0^{\,0:s-1}(\omega)\dd\omega & \dots & \displaystyle\int^h_0\ell^{\,0:s-1}_{s-1}(\omega)\dd\omega \end{pmatrix}
\]
	Now by \eqref{ab gp} and the standard formulae for Gaussian conditioning, we have
	\begin{align*}
	\mathbb{E}[y_{i+1}|y_i,& f_{i-s+1:i}] = \\
	&\phantom{=}
	\begin{pmatrix} \Phi_{i+1}^T\Phi_{i} \\ \Phi_{i+1}^T\phi_{i} \\ \vdots \\ \Phi_{i+1}^T\phi_{i-s+1} \end{pmatrix}^T
	\underbrace{\begin{pmatrix}
		\Phi_{i}^T\Phi_{i} & \Phi_{i}^T\phi_{i} & \cdots & \Phi_{i}^T\phi_{i-s+1} \\
		\phi_{i}^T\Phi_{i} & \phi_i^T\phi_i & \cdots & \phi_i^T\phi_{i-s+1}\\
		\vdots & \vdots & \ddots & \vdots \\
		\phi_{i-s+1}^T\Phi_{i} & \phi_{i-s+1}^T\phi_{i+1} & \dots & \phi_{i-s+1}^T\phi_{i-s+1}
	\end{pmatrix}^{-1}}_{\displaystyle\mathbb{I}_{s+1}^{-1}}
	\begin{pmatrix}  y_{i}\\ f_i\\ \vdots \\ f_{i-s+1} \end{pmatrix} \\
		&=  (\Phi_{i+1}^T\Phi_{i})y_i + \sum_{k=0}^{s-1}(\Phi_{i+1}^T\phi_{i-k})f_{i-k} \\
	& = y_i + \sum_{k=0}^{s-1}[\Phi_{i+1}]_{k+2}\cdot f_{i-k} \qquad \qquad \qquad \left( \parbox{12.5em}{where $[\Phi_{i+1}]_{k+2}$ denotes the \\\hspace{2em} $(k+2)$th component of $\Phi_{i+1}$} \right)\\
	&= y_i + \sum_{k=0}^{s-1} \left[ \int_0^h \ell^{\,0:s-1}_k(\omega)\dd\omega \right] \cdot f_{i-k} \\
	&= y_i + h\sum_{k=0}^{s-1} c_{k,s} f_{i-k} \qquad \qquad \qquad \qquad \quad \text{since} \int_0^h \ell^{\,0:s-1}_k(\omega)\dd\omega =  h c_{k,s}
	\end{align*}
which is equal to the $s$-step Adams-Bashforth predictor defined by \eqref{lagrange interpolation} and \eqref{ab def}. Next we write
\begin{align*}
	\renewcommand\arraystretch{1.2}
	\mathrm{Var}[y_{i+1}|y_i,f_{i-s+1:i}] &= \Phi_{i+1}^T\Phi_{i+1} -
	\begin{pmatrix} \Phi_{i+1}^T\Phi_{i} \\ \Phi_{i+1}^T\phi_{i} \\ \vdots \\ \Phi_{i+1}^T\phi_{i-s+1} \end{pmatrix}^T \mathbb{I}_{s+1}^{-1}
	\begin{pmatrix} \Phi_{i}^T\Phi_{i+1} \\ \phi_{i}^T\Phi_{i+1} \\ \vdots \\ \phi_{i-s+1}^T\Phi_{i+1} \end{pmatrix} \\
	&= \Phi_{i+1}^T\Phi_{i+1} - 
	\begin{pmatrix} 1 \\ [\Phi_{i+1}]_2 \\ \vdots \\ [\Phi_{i+1}]_{s+1} \end{pmatrix}^T \begin{pmatrix} 1 \\ [\Phi_{i+1}]_2 \\ \vdots \\ [\Phi_{i+1}]_{s+1} \end{pmatrix} \\
	&= \Phi_{i+1}^T\Phi_{i+1} - \Phi_{i+1}^T\Phi_{i+1} \\
	&= 0
	\end{align*}
	and the proposition follows.
\subsection{Proof of Proposition 2} \label{p2 proof}
	We follow the same reasoning as in Proposition \ref{ab det theorem}. Since the additional basis function at the end of $\phi^+_{i-k}$ is clearly zero at for all $0 \leq k \leq s-1$, each inner product of the form $\phi^{+T} \phi^+$, $\Phi^{+T} \phi^+$ and $\phi^{+T} \Phi^+$ is equal to the corresponding inner product $\phi^{T} \phi$, $\Phi^{T} \phi$ and $\phi^{T} \Phi$ as no additional contribution from the new extended basis arises. It therefore suffices to check only the terms of the form $\Phi^{+T}\Phi$. \\
	
	Integrating the additional basis function gives a polynomial of degree $s+1$ with a constant factor $\omega$. Evaluating this at $t_i = 0$ means that the additional term is also 0 in $\Phi_i$. Therefore $\Phi_{i+1}^{+T}\Phi_i^+ = \Phi_{i+1}^{T}\Phi_i$ and $\Phi_{i}^{+T}\Phi_i^+ = \Phi_{i}^{T}\Phi_i$. It follows that the expression for $\mathbb{E}[y_{i+1}|y_i,f_{i-s+1:i}]$ is exactly the same as when using the unaugmented basis function set.\\
	
	The argument in the previous paragraph means we can immediately write down that 
	\[
	\mathrm{Var}[y_{i+1}|y_i,f_{i-s+1:i}] = \Phi_{i+1}^{+T}\Phi_{i+1}^+ - \Phi_{i+1}^T\Phi_{i+1}
	\]
	Since the first $s+1$ components of $\Phi_{i+1}^{+T}$ are equal to the $s+1$ components of $\Phi_{i+1}^{T}$, this expression reduces to the contribution of the augmented basis element. Therefore
	\begin{align*}
	\mathrm{Var}[y_{i+1}|y_i,f_{i-s+1:i}] &= \left( \alpha h^{s} \displaystyle\int_0^h \ell^{\,-1:s-1}_{-1}(\omega) \dd\omega \right)^2 \\&=  \left( \alpha h^{s+1} \beta^{AM}_{-1,s+1} \right)^2
	\end{align*}
	The Adams-Moulton coefficient $\beta^{AM}_{-1,s+1}$ is equal to the local truncation error constant for the $s$-step Adams-Bashforth method \cite{butcher08} and the proposition follows.

\subsection{Extension to Adams-Moulton} \label{am extension}
We collect here the straightforward modifications required to the constructions in the main paper to produce implicit Adams-Moulton methods instead of explicit Adams-Bashforth versions.\\

The telescopic decomposition \eqref{product telescope} becomes 
\begin{equation}
	p(y_{1:N},f_{0:N}|y_0) = \prod_{i=0}^{N} p(f_{i}|y_{i}) \times \prod_{i=0}^{N-1} p (y_{i+1}|y_i,f_{\mathrm{max}(0,i-s+1):i+1})
\end{equation}
where it is particularly to be noted that $f_N$ is no longer superfluous.\\

The Lagrange interpolation resulting in the the Adams-Moulton method is
\begin{equation}
Q_i(\omega) = \sum_{j=-1}^{s-1} \ell_j^{\,-1:s-1}(\omega) f_{i-j}  \qquad\qquad \ell_j^{\,-1:s-1}(\omega) = \prod_{\substack{k=-1 \\ k \neq j}}^{s-1} \frac{\omega - t_{i-k}}{t_{i-j}-t_{i-k}} \label{lagrange AM},
\end{equation}
the analogous vectors of basis polynomials to \eqref{phi} and \eqref{Phi} are
\begin{align}
\psi(\omega) &= \begin{pmatrix} \vphantom{\displaystyle\int} 0 & \ell^{\,-1:s-1}_{-1}(\omega) & \ell^{\,-1:s-1}_0(\omega) & \ell^{\,-1:s-1}_1(\omega) & \dots & \ell^{\,-1:s-1}_{s-1}(\omega) \end{pmatrix}^T \label{psi}\\
\Psi(\omega) &= \int \psi(\omega) \dd \omega  = \begin{pmatrix} 1 & \displaystyle\int\ell^{\,-1:s-1}_{-1}(\omega)\dd\omega & \dots & \displaystyle\int\ell^{\,-1:s-1}_{s-1}(\omega)\dd\omega \end{pmatrix}^T \label{Psi} 
\end{align}
and the iterator is defined by
\begin{equation}
y_{i+1} - y_{i} \approx \int_{t_i}^{t_{i+1}} Q_i(\omega) \dd\omega = h \sum_{j=-1}^{s-1} \beta^{AM}_{j,s+1} f_{i-j} \label{am def}
\end{equation}
with $\beta^{AM}_{j,s+1} \equiv h^{-1}\displaystyle\int_0^h \ell_j^{\,-1:s-1}(\omega)\dd\omega$ are the Adams-Moulton coefficients.\\

The Gaussian process prior resulting in AM is

\begin{equation}
\renewcommand\arraystretch{1.2}
\begin{pmatrix}
	y_{i+1}\\ y_{i}\\ f_{i+1}\\ f_i\\ f_{i-1}\\ \vdots \\ f_{i-s+1}
\end{pmatrix} = \mathcal{N} \left[ \begin{pmatrix}
	0 \\ 0 \\ 0 \\ 0 \\0 \\ \vdots \\ 0 
\end{pmatrix}, \begin{pmatrix}
	\Psi_{i+1}^T\Psi_{i+1} & \Psi_{i+1}^T\Psi_{i} & \Psi_{i+1}^T\psi_{i+1} & \cdots & \Psi_{i+1}^T\psi_{i-s+1}\\
	
	\Psi_{i}^T\Psi_{i+1} & \Psi_{i}^T\Psi_{i} & \Psi_{i}^T\psi_{i+1} & \cdots & \Psi_{i}^T\psi_{i-s+1}\\
	
	\psi_{i+1}^T \Psi_{i+1} & \psi_{i+1}^T\Psi_i & \psi_{i+1}^T\psi_{i+1} & \cdots & \psi_{i+1}^T\psi_{i-s+1}\\
	
	\psi_{i}^T\Psi_{i+1} & \psi_{i}^T\Psi_{i} & \psi_{i}^T\psi_{i+1} & \cdots & \psi_{i}^T\psi_{i-s+1}
	 \\
	
	\psi_{i-1}^T\Psi_{i+1} & \psi_{i-1}^T\Psi_{i} & \psi_{i-1}^T\psi_{i+1} & \cdots & \psi_{i-1}^T\psi_{i-s+1}\\
	
	\vdots & \vdots & \vdots & \ddots & \vdots\\
	\psi_{i-s+1}^T\Psi_{i+1}&
	\psi_{i-s+1}^T\Psi_{i} & \psi_{i-s+1}^T\psi_{i+1} & \cdots & \psi_{i-s+1}^T\psi_{i-s+1}
\end{pmatrix} \right] \label{am gp}
\end{equation}

\subsection{Adams-Moulton integrator with $s=4$} \label{s4}
The conditional distribution of interest is $p(y_{i+1}|y_i,f_{i+1},f_i,f_{i-1},f_{i-2}) \equiv p(y_{i+1}|y_i,f_{i+1:i-2})$. In the deterministic case the vectors of basis functions become
\begin{align}
\psi(\omega)_{s=4} &= \begin{pmatrix} 0 & \frac{\omega(\omega+h)(\omega+2h)}{-6h^3} & \frac{(\omega - h)(\omega+h)(\omega+2h)}{2h^3} & \frac{\omega(\omega - h)(\omega+2h)}{-2h^3} & \frac{\omega(\omega - h)(\omega+h)}{6h^3} \end{pmatrix} \\
\Psi(\omega)_{s=4} &= \begin{pmatrix} 1 & \frac{\omega^2(2h + \omega)^2}{24h^3} & \frac{\omega(3\omega^3 + 8h\omega^2 - 6h^2\omega - 2fh^3)}{-24h^3} & \frac{\omega^2( 3\omega^2 + 4h\omega - 12h^2)}{24h^3} & \frac{\omega^2(\omega^2 - 2h^2)}{-24h^3} \end{pmatrix}
\end{align}
and the resulting calculations give 
\begin{align*}
\mathbb{E}(y_{i+1}|y_i,f_{i+1:i-2}) &= y_i + h\left( \dfrac{3}{8}f_{i+1} + \dfrac{19}{24}f_{i} -\dfrac{5}{24}f_{i-1} + \dfrac{1}{24}f_{i-2} \right) \\
\mathrm{\Var}(y_{i+1}|y_i,f_{i+1:i-2}) &= 0
\end{align*}
The probabilistic version is
\begin{align}
\psi^+(\omega)_{s=4} &= \begin{pmatrix}& \cdots \psi(\omega)_{s=4} \cdots & \frac{\alpha\omega(\omega-h)(\omega+h)(\omega+2h)}{24} \end{pmatrix} \\
\Psi^+(\omega)_{s=4} &= \begin{pmatrix}& \cdots \Psi(\omega)_{s=4} \cdots &  \frac{\alpha\omega^2 (6\omega^3 + 15\omega^2h - 10\omega h^2 - 30h^3)}{720} \end{pmatrix}
\end{align}
and further calculation shows that
\begin{align}
\mathbb{E}(y_{i+1}|y_i,f_{i-1:i+2}) &= y_i + h\left( \dfrac{3}{8}f_{i+1} + \dfrac{19}{24}f_{i} -\dfrac{5}{24}f_{i-1} + \dfrac{1}{24}f_{i-2} \right) \\
\mathrm{\Var}[y_{i+1}|y_i,f_{i+1:i-2}] &= \left(\dfrac{19 h^5 \alpha}{720} \right)^2 \label{variance of AB3}
\end{align}

\subsection*{Remark}
Proofs analogous to those of Propositions 1 and 2, for the Adams-Moulton case, follow the same line of reasoning as for the Adams-Bashforth case.

\subsection{Expansion of backward difference coefficient approximation for $s=3$} \label{bdc expansion}
From \eqref{bdc}, we have for $s=3$
\begin{align*}
y_{i+1} &= y_i + h\left( \dfrac{23}{12}f_i -\dfrac{4}{3}f_{i-1} + \dfrac{5}{12}f_{i-2} \right) - \dfrac{3}{8}h^4y^{(4)}(t_{i+1}) + O(h^{5}) \\
&= y_i + h\left( \dfrac{23}{12}f_i -\dfrac{4}{3}f_{i-1} + \dfrac{5}{12}f_{i-2} \right) - \dfrac{3}{8}h^4f'''(t_{i+1}) + O(h^{5}) \\
&= y_i + h\left( \dfrac{23}{12}f_i -\dfrac{4}{3}f_{i-1} + \dfrac{5}{12}f_{i-2} \right) - \dfrac{3}{8}h^4
\left[\dfrac{-f_i + 3 f_{i-1} - 3 f_{i-2} + f_{i-3}}{h^3} + O(h) \right] + O(h^{5})\\
&= \underbrace{y_i + h \left( \dfrac{55}{24}f_i -\dfrac{59}{24}f_{i-1} + \dfrac{37}{24}f_{i-2} - \dfrac{3}{8}f_{i-3} \right)}_{\mathrm{AB4}} + \;O(h^5)
\end{align*}

\subsection{Proof of Theorem 3} \label{t3 proof}

Proposition 2 implies that our integrator can be written as 
\begin{equation} \label{eq:num}
y_{i+1}=y_{i}+h\sum_{j=0}^{s-1}\beta^{AB}_{j,s}f(y_{i-j},t_{i-j})+ \xi_{i}
\end{equation}
where $y_i$ denotes the numerical solution at iteration $i$, and $\xi_{i} \in \mathbb{R}^d$ is a Gaussian random variable satisfying $\mathbb{E}|\xi_i \xi_i^T| = Qh^{2s+2}$ for some fixed $d \times d$ matrix $Q$. We denote the true solution of the ODE \eqref{ode defn} at iteration $i$ by $Y_{i} \equiv y(t_{i})$ and we have that 
\begin{equation} \label{eq:true}
Y_{i+1}=Y_{i}+h\sum_{j=0}^{s-1}\beta^{AB}_{j,s}f(Y_{i-j},t_{i-j})+ \tau_{i}
\end{equation}
where by construction the local truncation error $\tau_{i}=O(h^{s+1})$. If we now subtract \eqref{eq:num} from \eqref{eq:true} and denote the accumulated error at iteration $i$ by $E_{i}=Y_{i}-y_{i}$, we have 
\[
E_{i+1}=E_{i} +\Delta \phi_{i}+\tau_{i} - \xi_{i}
\]
where 
\[
\Delta \phi_{i} := h \sum_{j=0}^{s-1} \beta_{j,s}^{AB}\Delta f_{i-j}, \qquad  \Delta f_{i-j}:=f(Y_{i-j},t_{i-j})-f(y_{i-j},t_{i-j})
\]
We will rearrange this $s$-step recursion to give an equivalent one-step recursion in an higher-dimensional space. In particular, using the trivial identities 
$E_{i-1}=E_{i-1}, \cdots, E_{i-s+1}=E_{i-s+1}$ we obtain 
\[
\underbrace{\left( \begin{array}{c}
E_{i+1}  \\
E_{i}  \\
\vdots \\
E_{i-s+2}   
\end{array} \right)}_{\displaystyle =:\mathcal{E}_{i+1}}=\underbrace{\left( \begin{array}{cccc}
\mathbb{I}_d &  0 & \cdots & 0  \\
\mathbb{I}_d & 0 & \cdots & 0 \\
& \ddots & \ddots & \\
0 & &  \mathbb{I}_d & 0  
\end{array} \right)}_{\displaystyle =:\mathcal{A}}\underbrace{\left( \begin{array}{ccc}
E_{i}  \\
E_{i-1}  \\
\vdots \\
E_{i-s+1}   
\end{array} \right)}_{\displaystyle =: \mathcal{E}_{i}}+\underbrace{\left( \begin{array}{c}
\Delta \phi_{i}  \\
0  \\
\vdots \\
0   
\end{array} \right)}_{\displaystyle =:\Delta \Phi_{i}}+\underbrace{\left( \begin{array}{c}
\tau_{i}  \\
0  \\
\vdots \\
0   
\end{array} \right)}_{\displaystyle =:\mathcal{T}_{i}}-
\underbrace{\left( \begin{array}{c}
\xi_{i}  \\
0  \\
\vdots \\
0   
\end{array}\right)}_{\displaystyle =:\Xi_{i}}
\]
or in compact form,
\[
\mathcal{E}_{i+1}=\mathcal{A}\mathcal{E}_{i}+\Delta \Phi_{i}+\mathcal{T}_{i} - \Xi_i, \qquad i=s-1, \dots, N-1, \qquad N = T/h \numberthis \label{E_trace}
\]
For the subsequent calculations it will be necessary to find a scalar product inducing a matrix norm such that the norm of the matrix $\mathcal{A}$ 
is less or equal to $1$. This is possible if the eigenvalues  of the Frobenius  matrix $\mathcal{A}$ lie inside the unit circle on the complex plane 
and are simple if their modulus is equal to $1$. It is easy to show that the eigenvalues of $\mathcal{A}$ are roots of the characteristic polynomial associated with the deterministic integrator \eqref{ab def}. Since we have assumed that the  deterministic integrator is convergent, $\mathcal{A}$ does have the claimed property, since it is equivalent to the root condition in Dahlquist's equivalence theorem \cite{butcher08}. Thus there exists a non-singular matrix $\Lambda$ with a  block structure like $\mathcal{A}$ such that $||\Lambda^{-1} \mathcal{A} \, \Lambda ||_{2} \leq 1$. We can therefore choose a scalar product for $\mathcal{X},\mathcal{Y} \in \mathbb{R}^{ds}$ as
\[
\left \langle \mathcal{X},\mathcal{Y} \right \rangle_{*}:= \left \langle \Lambda ^{-1}\mathcal{X}, \Lambda^{-1}\mathcal{Y} \right \rangle_{2}
\]
and then have $|\cdot|_{*}$ and $|| \cdot||_{*}$ as the induced vector and matrix norms respectively, with $||\mathcal{A} ||_{*}=||\Lambda ^{-1}\mathcal{A} \,\Lambda ||_{2} \leq 1$ as required.  We also have 
\[
\left \langle \mathcal{X},\mathcal{Y} \right \rangle_{*}=\mathcal{X}^{T} \Lambda^{-T}\Lambda^{-1}\mathcal{Y}=\mathcal{X}^{T}\Lambda^{*}\mathcal{Y} \ \text{with} \ \Lambda^{*}=\Lambda^{-T}\Lambda^{-1}= (\lambda^{*}_{ij}\otimes\mathbb{I}_{d})_{1\leq i,j\leq s} \numberthis \label{lambda}
\]
Due to the equivalence of norms there exist constants $c^{*},c_{*}>0$ such that 
\[
|\mathcal{X}|^{2}_{2} \leq c^{*}|\mathcal{X}|^{2}_{*} \quad \text{and} \quad |\mathcal{X}|^{2}_{*} \leq c_{*}|\mathcal{X}|^{2}_{\infty} \quad \text{for all } \mathcal{X} \in \mathbb{R}^{ds}, \numberthis \label{norm equivalence}
\]
where 
$
|\mathcal{X}|^{2}_{2}=\sum_{j=1,\dots, s}|x_{j}|^{2} $ and $  |\mathcal{X}|_{\infty}=\max_{j=1,\dots,s}|x_{j}| 
$ for $\mathcal{X}=(x^{T}_{1},\cdots,x^{T}_{s})^{T}, \ x_j \in \mathbb{R}^d$.

For the particular vectors $\tilde{\mathcal{X}}=(x^{T},0,\cdots,0)^{T}$  and $\tilde{\mathcal{Y}}=(y^{T},0,\cdots,0)^{T}$ with $\tilde{\mathcal{X}},\tilde{\mathcal{Y}}\in \mathbb{R}^{ds}$ and $x,y \in \mathbb{R}^{d}$, one has
\[
\langle \tilde{\mathcal{X}},\tilde{\mathcal{Y}} \rangle_{*}=\lambda^{*}_{11}\langle x,y \rangle_{2}=\lambda^{*}_{11}x^{T}y, \numberthis \label{particular vectors}
\]
where $\lambda^{*}_{11}$ is as in \eqref{lambda}. Applying the norm $|\cdot|^{2}_{*}$ to \eqref{E_trace} and taking expectations gives
\begin{align*}
\mathbb{E}|\mathcal{E}_{i+1}|^2_* &= \mathbb{E}|\mathcal{A}\mathcal{E}_{i}+\Delta \Phi_{i}+\mathcal{T}_{i} - \Xi_i|^2_* \\
&= \mathbb{E}|\mathcal{A}\mathcal{E}_{i}+\Delta \Phi_{i}+\mathcal{T}_{i}|^2_* + O(h^{2s+2})\\
&= \mathbb{E}|\mathcal{A}\mathcal{E}_{i}+\Delta \Phi_{i}|^2_* +2 \mathbb{E}\langle h^{1/2}(\mathcal{A}\mathcal{E}_{i}+\Delta \Phi_{i}),\mathcal{T}_i h^{-1/2} \rangle_* + \mathbb{E}|\mathcal{T}_{i}|^2_* + O(h^{2s+2}) \\
&= \mathbb{E}|\mathcal{A}\mathcal{E}_{i}+\Delta \Phi_{i}|^2_* +2 \mathbb{E}\langle h^{1/2}(\mathcal{A}\mathcal{E}_{i}+\Delta \Phi_{i}),\mathcal{T}_i h^{-1/2} \rangle_* + O(h^{2s+2}) \numberthis \label{three terms}
\end{align*}
We now consider the term $|\mathcal{A}\mathcal{E}_{i}+\Delta \Phi_{i}|^2_*$ and expand it as
\begin{align*}
|\mathcal{A}\mathcal{E}_{i}+\Delta \Phi_{i}|^2_* &= 	\underbrace{|\mathcal{A}\mathcal{E}_{i}|^2_*}_{\text{A}} + \underbrace{|\Delta \Phi_{i} |_*^2}_{\text{B}} + \underbrace{2 \langle \mathcal{A}\mathcal{E}_{i},\Delta \Phi_{i} \rangle_*}_{\text{C}}
\end{align*}
For term A we immediately have $|\mathcal{A}\mathcal{E}_{i}|^2_* \leq |\mathcal{E}_{i}|^2_*$ by construction of the norm $|\cdot|^{2}_{*}$.

For term B we have that
\begin{align*}
	 |\Delta\Phi_i|^2_* &= \lambda_{11}^*|\Delta\phi_i|^2 && \text{from \eqref{particular vectors}}\\
	 &= \lambda_{11}^*\big| h \textstyle\sum_{j=0}^{s-1} \beta_{j,s}^{AB}\Delta f_{i-j} \big|^2 \\
	 &\leq \lambda_{11}^* s h^2 \textstyle\sum_{j=0}^{s-1} \left|  \beta_{j,s}^{AB}\Delta f_{i-j} \right|^2 &&\text{by Cauchy-Schwarz}\\
	 &\leq \lambda_{11}^* s h^2 L_f^2 \textstyle\sum_{j=0}^{s-1}  (\beta_{j,s}^{AB})^2 |E_{i-j}|^2 &&\text{since $f$ is Lipschitz} \\
	 &\leq \lambda_{11}^* s h^2 L_f^2 C_\beta^2 \textstyle\sum_{j=0}^{s-1} |E_{i-j}|^2 &&\text{where } C_\beta^2 = \max_{j=0,\dots,s-1}\beta_{j,s}^{AB} \\
	 &\leq \lambda_{11}^* s h^2 L_f^2 C_\beta^2 c^* |\mathcal{E}_{i}|^2_* &&\text{from \eqref{norm equivalence}} \\
	 &=: \Gamma^2 h^2 |\mathcal{E}_{i}|^2_* && \text{where } \Gamma^2 = \lambda_{11}^* s L_f^2 C_\beta^2 c^*
\end{align*}
 
For term C we have $ 2 \langle \mathcal{A}\mathcal{E}_{i},\Delta \Phi_{i} \rangle_* \leq 2 |\mathcal{A}\mathcal{E}_{i}|_* |\Delta\Phi_i|_* \leq 2\Gamma h |\mathcal{E}_{i}|^2_* $ and it follows that
\[
|\mathcal{A}\mathcal{E}_{i}+\Delta \Phi_{i}|^2_* \leq (1+O(h)) |\mathcal{E}_{i}|^2_*
\]
Then from \eqref{three terms} we have
\begin{align*}
\mathbb{E}|\mathcal{E}_{i+1}|^2_* &= 	\mathbb{E}|\mathcal{A}\mathcal{E}_{i}+\Delta \Phi_{i}|^2_* +2 \mathbb{E}\langle h^{1/2}(\mathcal{A}\mathcal{E}_{i}+\Delta \Phi_{i}),\mathcal{T}_i h^{-1/2} \rangle_* + O(h^{2s+2}) \\
& \leq (1+O(h)) \mathbb{E} |\mathcal{E}_{i}|^2_* + 2h \mathbb{E}|\mathcal{A}\mathcal{E}_{i}+\Delta \Phi_{i}|^2_*  + 2h^{-1}\mathbb{E} |\mathcal{T}_i|^2_* +  O(h^{2s+2}) \\
& \leq (1+O(h)) \mathbb{E} |\mathcal{E}_{i}|^2_* +  O(h^{2s+1}) + O(h^{2s+2}) \numberthis \label{pre-gronwall}
\end{align*}
Then by applying the Gronwall inequality we have (for different $K$ in each line)
\[
\max_{0\leq kh \leq T} \mathbb{E}|\mathcal{E}_k|^2_* \leq K(T)h^{2s}
\]
and since $\mathcal{E}_k = (E_k,E_{k-1},\cdots,E_{k-s+1})$ we conclude that
\[
\max_{0\leq kh \leq T} \mathbb{E}|E_k|^2 \leq K(T) h^{2s}
\]
Note that in \eqref{pre-gronwall}, the $O(h^{2s+2})$ term derived from the introduced perturbations $\xi_i$ is of one higher order than the $O(h^{2s+1})$ term representing the truncation error in the deterministic solver. This observation implies that a noise vector satisfying $\mathbb{E}|\xi_i \xi_i^T| = Qh^{2s+1}$ would also give rise to an integrator of order $s$.
\subsection*{Remark}
An analogous proof for the Adams-Moulton case follows with straightforward modifications.

\end{document}